\documentclass[11pt]{article}
\usepackage{amsmath}
\usepackage{amssymb}
\usepackage{theorem}
 \newtheorem{thm}{Theorem}[section]

 \theoremstyle{definition}
 
 \theoremstyle{remark}
 \newtheorem{rem}[thm]{Remark}
 
 \numberwithin{equation}{section}
\usepackage{mathtools}
\usepackage{synttree}
\usepackage{subfigure}

\usepackage{xcolor}

\newcommand{\green}[1]{#1}

\newcommand{\Rbb}{\mathbb{R}}

\newcommand{\Nbb}{\mathbb{N}}
\newcommand{\Rc}{\mathcal{R}}

\newcommand{\Mc}{\mathcal{M}}

\newcommand{\Dc}{\mathcal{D}}

\newcommand{\Jc}{\mathcal{J}}

\newcommand{\Kc}{\mathcal{K}}

\newcommand{\rank}{\mathrm{rank}}

\begin{document}

\title{\bfseries Low-rank tensor methods for model order reduction}

\author{Anthony Nouy\footnotemark[2]}

\renewcommand{\thefootnote}{\fnsymbol{footnote}}
\footnotetext[2]{Ecole Centrale de Nantes, GeM, UMR CNRS 6183, France. e-mail: \texttt{anthony.nouy@ec-nantes.fr}
}

\renewcommand{\thefootnote}{\arabic{footnote}}

\date{}

\maketitle

%



\begin{abstract}
Parameter-dependent models arise in many contexts such as uncertainty quantification, sensitivity analysis, inverse problems or optimization. Parametric or uncertainty analyses usually require the evaluation of an output of a model for many instances of the input parameters, which may be intractable for complex numerical models. A possible remedy consists in replacing the model by an approximate model with reduced complexity (a so called reduced order model) allowing a fast evaluation of output variables of interest. 
This chapter provides an overview of low-rank methods for the approximation of functions that are identified  either with order-two tensors (for vector-valued functions) or higher-order tensors (for multivariate functions).
Different approaches are presented for the computation of low-rank approximations, either based on samples of the function or on the equations that are satisfied by the function, the latter approaches including projection-based model order reduction methods. For multivariate functions, different notions of ranks and the corresponding low-rank approximation formats are introduced.
\bigskip

\noindent \textbf{Keywords}: uncertainty quantification, high-dimensional problems, model order reduction, low-rank approximation, parameter-dependent equations, stochastic equations. 
\end{abstract}

\maketitle

\newpage

\section*{Introduction}

Parameter-dependent models arise in many contexts such as uncertainty quantification, sensitivity analysis, inverse problems and optimization. 
These models are typically given under the form of a \green{parameter-dependent equation}
\begin{equation}
R(u(\xi);\xi) = 0, \label{residualequation}
\end{equation}
where $\xi$ are parameters taking values in some set $\Xi$, 
and where the solution 
 $u(\xi)$ is in some vector space $V$, say $\Rbb^M$. Parametric or uncertainty analyses usually require the evaluation of the solution for many instances of the parameters, which may be intractable for complex numerical models (with large $M$) for which one single solution requires hours or days of computation time. Therefore, one usually relies on  approximations of the solution map $u:\Xi \to V$ allowing for a rapid evaluation of output quantities of interest. These approximations take on different names such as 
meta-model, surrogate model or reduced order model.
They are usually of the form 
\begin{equation}
u_m(\xi) = \sum_{i=1}^m v_i s_i(\xi),\label{approxrankm}
\end{equation}
where the $v_i $ are elements in $V$ and the $s_i$ are elements of some space $S$ of functions defined on $\Xi$. Standard linear approximation methods rely on the introduction of generic bases (e.g. polynomials, wavelets, etc) allowing an accurate approximation of a large class of models to be constructed but at the price of requiring expansions with a high number of terms $m$. These generic approaches usually result in a very high computational complexity.
\\\par
\green{Model order reduction}  methods  aim at finding an approximation $u_m$ with a small number of terms ($m\ll M$)
 that are adapted to the particular function $u$. 
One  can distinguish approaches relying (i) on the construction of a reduced basis $\{v_1,\hdots,v_m\}$ in $V$, (ii) on the construction of a reduced basis $\{s_1,\hdots,s_m\}$ in $S$ or (iii) directly on the construction of  an approximation under the form \eqref{approxrankm}. These approaches are closely related. They all result in a \green{\emph{low-rank approximation}} $u_m$, which can be interpreted as a rank-$m$ element of the tensor space $V\otimes S$. 
Approaches of type (i) are usually named \green{\emph{Projection-based model order reduction methods}} since they define $u_m(\xi)$ as a certain projection 
of $u(\xi)$ onto a low-dimensional subspace of $V$. 
They include \green{\emph{Reduced Basis}}, \green{\emph{Proper Orthognonal Decomposition}}, \emph{Krylov subspace}, \emph{Balanced truncation methods}, and also subspace-based variants of \green{\emph{Proper Generalized Decomposition}} methods. Corresponding reduced order models usually take the form of a small system of equations which defines 
the projection $u_m(\xi)$ for each instance of $\xi$. Approaches (i) and (ii) are in some sense dual to each other. Approaches (ii) include
\green{\emph{sparse approximation methods}} which consist in selecting $\{s_1,\hdots,s_m\}$ in a certain dictionary of functions (e.g. a polynomial  basis), based on prior information on $u$ or based on a posteriori error estimates (adaptive selection). They also include methods for the construction of reduced bases in $S$ that exploit some prior information on $u$.
\emph{Low-rank tensor methods} enter the family of approaches (iii), where approximations of the form \eqref{approxrankm} directly result from an optimization on \green{low-rank manifolds}. \par
When one is interested not only in evaluating the solution $u(\xi)$ at a finite number of samples of $\xi$ but in obtaining an explicit representation of the solution map $u:\Xi \to V$, approximations of functions of multiple parameters $\xi=(\xi_1,\hdots,\xi_d)$ are required. This constitutes a challenging issue for high-dimensional problems.  Naive approximation methods which consist in using tensorized bases 
yield  an exponential increase in storage or computational complexity when the dimension $d$ increases, which is the so-called 
\green{\emph{curse of dimensionality}}. Specific structures of the functions have to be exploited in order to reduce the complexity. Standard structured approximations include additive approximations
$ u_1(\xi_1) + \hdots + u_d(\xi_d)$, separated approximations $
u_1(\xi_1) \hdots u_d(\xi_d),$
 sparse approximations, or rank-structured approximations. Rank-structured approximation methods include several types of approximation depending on the notion of rank.  
\\
\par
The present chapter provides an overview of model order reduction methods based on low-rank tensor approximation methods.
In a first section, we recall some basic notions about low-rank approximations of an order-two tensor $u\in V\otimes S$ are recalled.  In the second section, which is devoted to  projection-based model reduction methods, we present different definitions of projections 
onto subspaces and different possible constructions of these subspaces. 
In the third section, we introduce the basic concepts of low-rank tensor methods for the approximation of a multivariate function, which is identified with a high-order tensor. The last two sections present different methods for the computation of low-rank approximations, either based on samples of the function (fourth section) or on the equations satisfied by the function (fifth section). 

\section{Low-rank approximation of order-two tensors}

Let us assume that $u: \Xi \to V$ with $V$ a Hilbert space equipped with a norm $\Vert \cdot\Vert_V$ and inner product $(\cdot,\cdot)_V$. For the sake of simplicity, let us consider that $V=\Rbb^M$. Let us further assume that $u$ is in the Bochner space $L^p_\mu(\Xi;V)$ for some $p\ge 1$, where $\mu$ is a probability measure supported on $\Xi$. 
The space $L^p_\mu(\Xi;V)$ can be identified with 
the algebraic tensor space $V\otimes L^p_\mu(\Xi)$ (or the completion of this space when $V$ is infinite dimensional, see e.g. \cite{Defant:1993fk})
which is the space of functions $w$ that can be written under the form 
$
w(\xi)   = \sum_{i=1}^m v_i  s_i (\xi)$ 
 for
 some $v_i\in V$ and $s_i\in L^p_\mu(\Xi)$, and some $m\in \Nbb$. 
 The rank of an element $w$, denoted $\rank(w)$, is the minimal integer $m$ such that $w$ admits such an $m$-term representation.
  A rank-$m$ approximation of $u$  then takes the form 
 \begin{equation}
u_m(\xi) = \sum_{i=1}^m v_i s_i(\xi),\label{eq:approx_um}
\end{equation}
and can be interpreted as a certain projection of $u(\xi)$ onto an $m$-dimensional subspace $V_m$ in $V$, where $\{v_1,\hdots,v_m\}$ constitutes a basis of $V_m$.

 \subsection{Best rank-$m$ approximation and optimal subspaces} 
The set of elements in  $V\otimes L^p_\mu(\Xi)$ with a rank bounded by $m$ is denoted 
 $  \Rc_m = \{w \in V\otimes L^p_\mu(\Xi) : \rank(w)\le m\}.
  $
The definition of a best approximation of $u$ from $\Rc_m$ requires the introduction of a measure of error. The best rank-$m$ approximation with respect to the  Bochner norm
$\Vert \cdot\Vert_p$ in  $L^p_\mu(\Xi;V)$ is the solution of 
  \begin{align}
  \min_{v \in \Rc_m} \Vert u - v \Vert_p = \min_{v\in \Rc_m} \Vert \Vert u(\xi) - v(\xi) \Vert_V \Vert_{L^p_\mu(\Xi)} := d_m^{(p)}(\Xi).\label{eq:best-rank-m}
  \end{align}
The set $\Rc_m$ admits the subspace-based parametrization
$\Rc_m = 
 \{w \in V_m\otimes L^p_\mu(\Xi) : V_m \subset V , \dim(V_m)=m\}.
$ 
 Then the best rank-$m$ approximation problem can be reformulated as an optimization problem 
over the set of $m$-dimensional spaces:
\begin{equation}
d_m^{(p)}(u)= \min_{\dim(V_m) = m } \min_{v\in V_m\otimes L^p_\mu} \Vert u-v \Vert_p = 
\min_{\dim(V_m)=m} \Vert u -P_{V_m} u\Vert_p  , \label{eq:optimal-space-norm}
\end{equation}
where $P_{V_m}$ is the orthogonal projection from $V$ to $V_m$ which provides the best approximation $P_{V_m} u(\xi)$ of $u(\xi)$ from $V_m$, defined by 
\begin{equation}
\Vert u(\xi)-P_{V_m}u(\xi) \Vert_V =  \min_{ v \in V_m}  \Vert u(\xi) - v \Vert_V
.\label{eq:PVm}
\end{equation}
That means that the best rank-$m$ approximation problem is equivalent to the problem of finding an optimal subspace of dimension $m$ for the projection of the solution. 

\subsection{Characterization of optimal subspaces}
The numbers $d_m^{(p)}(u)$, which are the best rank-$m$ approximation errors with respect to norms $\Vert \cdot\Vert_p$, are so-called \emph{linear widths} of the 
 \emph{solution manifold} 
$$\Kc := u(\Xi) = \{u(\xi) : \xi\in \Xi\} ,$$
and measure  how well the set of solutions $\Kc$ can be approximated by $m$-dimensional subspaces. They provide a quantification of the ideal performance of model order reduction methods. 

For $p=\infty$, assuming that $\Kc$ is compact, the number
\begin{equation}
d_m^{(\infty)}(u) = \min_{\dim(V_m)=m} \sup_{\xi\in \Xi} \Vert u(\xi) - P_{V_m}u(\xi) \Vert_V = \min_{\dim(V_m)=m} \sup_{v \in \Kc} \Vert v - P_{V_m} v \Vert_V \label{dminfini}\end{equation}
corresponds to the \emph{Kolmogorov $m$-width} $d_m(\Kc)_V$ of $\Kc$. This measure of error is particularly pertinent if one is interested in computing approximations that are uniformly accurate over the whole parameter set $\Xi$. 
For $p<\infty$, the numbers
\begin{align}
d_m^{(p)}(u) = \min_{\dim(V_m)=m} \left( \int_{\Xi} \Vert u(\xi) - P_{V_m}u(\xi) \Vert_V^p \mu(d\xi)  \right)^{1/p}\label{dmp}
\end{align}
provide measures of the error  that take into account the measure $\mu$ on the parameter set. This is particularly relevant for uncertainty quantification, where these 
numbers $d_m^{(p)}(u)$ directly control the error on the moments of the solution map $u$ (such as mean, variance or higher order moments). 
  
  Some general results on the convergence of linear widths are available in approximation theory (see e.g. \cite{Pietsch:1987_eigenvalues}). 
More interesting results which are specific to some classes of parameter-dependent equations have been recently obtained \cite{Lassila:2013uq,2015arXiv150206795C}. These results usually exploit smoothness and anisotropy of the solution map $u:\Xi \to V$ and are typically upper bounds deduced from results on polynomial approximation. 
Even if a priori estimates for  $d_m^{(p)}(u)$ are usually not available, a challenging problem is to propose numerical methods that provide approximations $u_m$ of the form \eqref{eq:approx_um} with an error  $\Vert u-u_m\Vert_p$ of the order of 
the best achievable accuracy $d_m^{(p)}(u)$.
 
 \subsection{Singular value decomposition}
Of particular importance is the case $p=2$ where $u$ is in the Hilbert space $L^2_\mu(\Xi;V)=V\otimes L^2_\mu(\Xi)$ and can be identified 
 with a compact operator $U : v \in V \mapsto  (u(\cdot),v)_V \in  L^2_\mu(\Xi)$ which admits 
a \green{\emph{singular value decomposition}} $$U = \sum_{i\ge 1} \sigma_i v_i  \otimes s_i,$$ where the $\sigma_i$ are the singular values and where $v_i$ and $s_i$ are the corresponding
normalized right and left singular vectors respectively. Denoting by $U^*$ the adjoint operator of $U$, defined by $U^* : s\in L^2_\mu(\Xi) \mapsto \int_\Xi  u(\xi) s(\xi) d\mu(\xi) \in V$, the operator 
 \begin{equation}
 U^*U := C(u) : v \in V\mapsto \int_\Xi u(\xi) (u(\xi),v)_V d\mu(\xi) \in V \label{eq:correlation_operator}
 \end{equation}
  is the correlation operator of $u$, with eigenvectors $v_i$ and corresponding eigenvalues $\sigma_i^2$.
  Assuming that singular values are sorted in decreasing order, $V_m = \mathrm{span}\{v_1,\hdots,v_m\}$ is a solution of 
 \eqref{eq:best-rank-m} (which means an optimal subspace) and a best approximation of the form \eqref{approxrankm}  is given by a rank-$m$  truncated singular value decomposition 
$u_m(\xi)=\sum_{i=1}^m \sigma_i v_i   s_i(\xi)$ which satisfies
 $
\Vert u - u_m \Vert_2 = d_m^{(2)}(u) = \big(\sum_{i>m} \sigma_i^2\big)^{1/2}.
 $ 
 
\section{Projection-based model order reduction methods}
Here we adopt a subspace point of view for the low-rank approximation of the solution map $u:\Xi\to V$. 
We first describe how to define projections onto a given subspace. 
Then we present methods for the practical construction of subspaces.

\subsection{Projections on a given subspace}
Here, we consider that a finite-dimensional subspace $V_m$ is given to us.
The best approximation of $u(\xi)$ from $V_m$ is given by the projection 
$P_{V_m} u(\xi)$ defined by \eqref{eq:PVm}.
However, in practice, projections of the solution $u(\xi)$ onto a given subspace $V_m$ must be defined using computable information.

\subsubsection{Interpolation}
When the subspace $V_m$ is the linear span of evaluations (samples) of the solution $u$ at $m$ given points  $\{\xi^1,\hdots,\xi^m\}$ in $ \Xi$, i.e.
 \begin{equation}
V_m = \mathrm{span}\{u(\xi^1),\hdots,u(\xi^m)\},
\label{eq:Vm_samples}
\end{equation} 
projections $u_m(\xi)$ of $u(\xi)$ onto $V_m$ can be obtained by interpolation. 
An interpolation 
of $u$ can be written in the form 
$$
u_m(\xi) =  \sum_{i=1}^m u(\xi^i) s_i(\xi),
$$ 
where functions $s_i$ satisfy the interpolation conditions   $s_i(\xi^j)=\delta_{ij}$, $1\le i,j\le m$. The best approximation $P_{V_m}u(\xi)$ is a particular case of interpolation. However, its computation is not of practical interest since it requires the knowledge of $u(\xi)$.
Standard polynomial interpolations can be used when $\Xi$ is an interval in $\Rbb$. In higher dimensions, polynomial interpolation can still be constructed for structured interpolation grids. For arbitrary sets of points, other interpolation formulae can be used, such as Kriging, nearest neighbor, Shepard or Radial Basis interpolations. 
These standard interpolation formulae provide approximations $u_m(\xi)$ that only depend on the value of $u$ at points $\{\xi^1,\hdots,\xi^m\}$.

Generalized interpolation formulae that take into account the function over the whole parameter set can be defined by
\begin{equation}
u_m(\xi) =  \sum_{i=1}^m u(\xi^i)  \varphi_i(u(\xi)),\label{geim}
\end{equation}
 with $\{\varphi_i\}$ a dual system to $\{u(\xi^i)\}$ such that $ \varphi_i(u(\xi^j ))=\delta_{ij}$ for  $1\le i,j\le m$.
This yields an interpolation $u_m(\xi)$ depending not only on the value of $u$ at the points   $\{\xi^1,\hdots,\xi^m\}$ but also on $u(\xi)$. This is of practical interest if 
$ \varphi_i(u(\xi))$ can be efficiently computed without a complete knowledge of $u(\xi)$. For example, if the coefficients of $u(\xi)$ on some basis of $V$ can be estimated without computing $u(\xi)$, then a possible choice consists in taking for $\{\varphi_1,\hdots,\varphi_m\}$ 
a set of functions that associate to an element of $V$ a set of $m$ of its coefficients. This is the idea behind the \green{\emph{Empirical Interpolation Method}} \cite{MAD09} and its generalization \cite{Maday:2013kx}.

\subsubsection{Galerkin projections}
For models described by an equation of type \eqref{residualequation}, the most prominent methods are \green{\emph{Galerkin projections}} which define the approximation $u_m(\xi)$ from the residual $R(u_m(\xi);\xi)$, e.g. by imposing orthogonality of the residual with respect to an $m$-dimensional space or by minimizing some  residual norm.
Galerkin projections do not provide the best approximation but under some usual assumptions, they can provide  quasi-best approximations $u_m(\xi)$ satisfying 
\begin{equation}
\Vert u(\xi)-u_m(\xi) \Vert_V \le c(\xi) \min_{v \in V_m} \Vert u(\xi) - v \Vert_V, \label{quasioptimality}
\end{equation}
where $c(\xi)\ge 1$. As an example, let us consider the case of a linear problem where 
\begin{equation}
R(v;\xi) = A(\xi) v - b(\xi),\label{residual_linear}
\end{equation} with $A(\xi)$ a linear operator from $V$ into some Hilbert space $W$, and let us consider a Galerkin projection defined by 
minimizing some residual norm, that means 
\begin{equation}
\Vert A(\xi) u_m(\xi) - b(\xi) \Vert_W = \min_{v \in V_m} \Vert A(\xi) v - b(\xi) \Vert_W. \label{eq:galerkin_minres}
\end{equation}
Assuming \begin{equation}
\alpha(\xi) \Vert v \Vert_V \le \Vert A(\xi) v \Vert_W \le \beta(\xi) \Vert v \Vert_V,\label{assumption_A}
\end{equation} the resulting projection $u_m(\xi)$ 
satisfies  \eqref{quasioptimality} with a constant $c(\xi) = \frac{\beta(\xi)}{\alpha(\xi)}$ which can be interpreted as the condition number of the operator $A(\xi)$. This 
reveals the interest of introducing efficient preconditioners in order to better exploit the approximation power of the subspace $V_m$ (see e.g. \cite{2014arXiv1403.7273C,2015arXiv150407903Z} for the construction of parameter-dependent preconditioners). The resulting approximation can be written $$u_m(\xi) = P_{V_m}^G(\xi) u(\xi),$$
where $P_{V_m}^G(\xi)$ is a parameter-dependent projection from $V$ onto $V_m$. Note that if $V_m$ is generated by samples of the solution, as in \eqref{eq:Vm_samples}, then the Galerkin projection is also an  interpolation which can be written under the form \eqref{geim} with parameter-dependent functions $\varphi_i$ that depend on the operator.
 
Some technical assumptions on the residual are required for a practical computation of Galerkin projections. More precisely, for $u_m(\xi) = \sum_{i=1}^m v_i s_i(\xi)$, the residual should admit a low-rank representation
$$
R(u_m;\xi) = \sum_{j} R_j \gamma_j(s(\xi);\xi),
$$
where the $R_j$ are independent of $s(\xi) = (s_1(\xi),\hdots,s_m(\xi))$ and where 
the $\gamma_j$ can be computed with a complexity depending on $m$ (the dimension of the reduced order model) but not on  the dimension of $V$. This allows Galerkin projections to be computed by solving a reduced system of $m$ equations on the unknown coefficients $s(\xi)$, 
with a computational complexity independent of the dimension of the full order model. 
For linear problems, such a property is obtained when the operator $A(\xi)$ and the right-hand side $b(\xi)$ admit low-rank representations (so-called \emph{affine representations})
\begin{equation}
A(\xi) = \sum_{i=1}^L {A_i} \alpha_i(\xi), \quad b(\xi)=  \sum_{i=1}^R b_i  \beta_i(\xi).\label{eq:affine_A_b}
\end{equation}
If $A(\xi)$ and  $b(\xi)$ are not explicitly given under the form \eqref{eq:affine_A_b}, then 
a preliminary approximation step is needed. Such expressions can   be obtained by the Empirical Interpolation Method
\cite{BAR02,Casenave:2014fk} or other low-rank truncation methods. Note that preconditioners for parameter-dependent operators should also have such representations in order to preserve a reduced order model with a complexity independent of the dimension of the full order model. 

 \subsection{Construction of subspaces}

The computation of an optimal $m$-dimensional subspace $V_m$ with respect to the natural norm in $L^p_\mu(\Xi;V)$, defined by \eqref{eq:optimal-space-norm}, is not feasible in practice since it requires the knowledge of $u$. Practical constructions of subspaces must rely on computable information on $u$, which can be samples of the solution or the model equations (when available).

\subsubsection{From samples of the function}
Here, we present constructions of subspaces $V_m$ which are based on evaluations of the function $u$ at a 
set of points $\Xi_K=\{\xi^1,\hdots,\xi^K\}$ in $\Xi$. Of course, $V_m$ can be chosen as the span of evaluations of $u$ at $m$ points chosen independently of $u$, e.g. through random sampling. Here, we present methods for obtaining subspaces $V_m$ that are closer to the optimal $m$-dimensional spaces, with $K\ge m$.

\paragraph{$L^2$ optimality.}
When one is interested in optimal subspaces with respect to the norm $\Vert\cdot\Vert_2$,  
the definition \eqref{eq:optimal-space-norm} of optimal subspaces can be replaced 
by \begin{equation}
\min_{\dim(V_m)=m}   \frac{1}{K} \sum_{k=1}^K    \Vert u(\xi^k) - P_{V_m} u(\xi^k) \Vert_{V}^2 
= \min_{v \in \Rc_m}   \frac{1}{K} \sum_{k=1}^K  \Vert u(\xi^k) - v(\xi^k) \Vert_{V}^2 ,\label{VmsamplesL2}
\end{equation}
where $\xi^1,\hdots,\xi^K$ are $K$ independent random samples drawn according the probability measure $\mu$. The resulting subspace $V_m$ is the dominant 
eigenspace of the empirical correlation operator 
\begin{equation}
C_K(u) = \frac{1}{K} \sum_{k = 1}^K u(\xi^k) (u(\xi^k),\cdot)_V,\label{estimatecorrelation}
\end{equation}
which is a statistical estimate of the  correlation operator $C(u)$ defined by \eqref{eq:correlation_operator}. 
This approach, which requires the computation of $K$ evaluations of the solution $u$, corresponds to the classical \green{\emph{Principal Component Analysis}}. It is at the basis of \green{Proper Orthogonal Decomposition} methods for parameter-dependent equations \cite{kahlbacher2006galerkin}. 
A straightforward generalization consists in defining the subspace $V_m$ by
  \begin{equation}
\min_{\dim(V_m)=m}  \sum_{k=1}^K \omega_k  \Vert u(\xi^k) - P_{V_m} u(\xi^k) \Vert_{V}^2 
= \min_{v \in \Rc_m}   \sum_{k=1}^K  \omega_k \Vert u(\xi^k) - v(\xi^k) \Vert_{V}^2 ,
\end{equation}
using a suitable quadrature rule for the integration over $\Xi$, e.g. exploiting the smoothness of the solution map $u : \Xi \to V$ in order to improve the convergence with $K$ and therefore decrease the number of evaluations of $u$ for a given accuracy. The resulting subspace $V_m$ is obtained as the dominant eigenspace of the operator 
\begin{equation}
C_K(u) =  \sum_{k = 1}^K \omega^k u(\xi^k) (u(\xi^k),\cdot)_V.\label{estimatecorrelationgeneralized}
\end{equation}

\paragraph{$L^\infty$ optimality.}
If one is interested in optimality with respect to the norm $\Vert\cdot\Vert_\infty$, the definition \eqref{eq:optimal-space-norm} of optimal spaces 
can be replaced by
\begin{equation}
\min_{\dim(V_m)=m}   \sup_{\xi \in \Xi_K}   \Vert u(\xi) - P_{V_m} u(\xi) \Vert_{V}
= \min_{v \in \Rc_m}    \sup_{\xi \in \Xi_K} \Vert u(\xi) - v(\xi) \Vert_{V}, \label{optimspaceLinf}
\end{equation}
with $\Xi_K=\{\xi^1,\hdots,\xi^K\}$ a set of $K$ points in $\Xi$.
A computationally tractable definition of subspaces can be obtained by adding  the constraint that subspaces $V_m$ are 
generated from $m$ samples of the solution, that means $V_m = \mathrm{span}\{u(\xi^1_*),\hdots,u(\xi^m_*)\}$. Therefore, problem \eqref{optimspaceLinf} becomes 
$$
\min_{\xi^1_*,\hdots,\xi^m_* \in \Xi_K} \max_{\xi\in \Xi_K} \Vert u(\xi) - P_{V_m} u(\xi) \Vert_V,
$$
where the $m$ points are selected in the finite set of points $\Xi_K$.  In practice, this combinatorial problem can be replaced by a \green{\emph{greedy algorithm}}, which consists in selecting the points  adaptively: given the first $m$ points and the corresponding subspace $V_m$, a new interpolation point $\xi^{m+1}_*$ is defined such that
\begin{equation}
 \Vert u(\xi^{m+1}_*) - P_{V_m} u(\xi^{m+1}_*) \Vert_V =\max_{\xi\in \Xi_K} \Vert u(\xi) - P_{V_m} u(\xi) \Vert_V.\label{eim}
\end{equation} 
This algorithm corresponds to the {Empirical Interpolation Method} \cite{BAR02,MAD09,Bebendorf:2014uq}.

\subsubsection{From approximations of the correlation operator}
An optimal $m$-dimensional space $V_m$ for the approximation of $u$ in $L^2_\mu(\Xi;V)$ is given by a dominant eigenspace of the 
correlation operator $C(u)$ of $u$. Approximations of optimal subspaces can then be obtained by computing dominant  eigenspaces of an approximate correlation operator.
These approximations can be obtained by using numerical integration in the definition of the correlation operator, yielding  
approximations \eqref{estimatecorrelation} or \eqref{estimatecorrelationgeneralized} which require evaluations of $u$. Another approach consists in 
using the correlation operator $C(u_K)$ (or the singular value decomposition) of an approximation $u_K$ of $u$ which can be obtained by a projection of $u$ onto a low-dimensional subspace $V\otimes S_K$ in $V\otimes L^2_\mu(\Xi)$. For example, an approximation can be sought after in the form $ u_K(\xi) = \sum_{k=1}^K v_k \psi_k(\xi)$ where   $\{\psi_k(\xi)\}_{k=1}^K$ is a polynomial basis. Let us note that the statistical estimate $C_K(u) $ in \eqref{estimatecorrelation} can be interpreted as the correlation operator $C(u_K)$ 
of a piecewise constant interpolation $u_K(\xi) = \sum_{k=1}^K u(\xi^k) 1_{\xi \in O_k}$ of $u(\xi)$, where 
the sets $\{O_1,\hdots,O_K\}$ form a partition of $\Xi$ such that $\xi^k \in O_k$ and $\mu(O_k) = 1/K$ for all $k$.

 \begin{rem}
Let us mention that 
the optimal rank-$m$ singular value decomposition of $u$ can be equivalently  obtained by computing the dominant eigenspace of the operator $\hat C(u) = UU^* : L^2_\mu(\Xi) \to L^2_\mu(\Xi)$. Then, a dual approach for model order reduction  consists in defining a subspace $S_m$ in $L^2_\mu(\Xi)$ as the dominant eigenspace of $\hat C( u_K)$, where $u_K$ is an approximation of $u$. An approximation of $u$ can then be obtained by a Galerkin projection onto the subspace $V\otimes S_m$ (see e.g. \cite{DOO07} where an approximation $u_K$ of the solution of a parameter-dependent partial differential equation is first computed using a coarse finite element approximation).
\end{rem}

\subsubsection{From the model equations} 
In the definition \eqref{eq:optimal-space-norm} of optimal spaces, $\Vert u(\xi) - v(\xi) \Vert$ can be replaced by a function $\Delta(v(\xi);\xi)$ which is computable without having $u(\xi)$ and such that $v\mapsto \Delta(v;\xi)$ has $u(\xi)$ as a minimizer over $V$. The choice for $\Delta$ is natural for problems where $u(\xi)$ is the minimizer of a functional 
$\Delta(\cdot;\xi) : V\to \Rbb$.  When $u(\xi)$ is solution of an equation of the form \eqref{residualequation}, 
a typical choice consists in taking for $\Delta(v(\xi);\xi) $ a certain norm  of the residual $R(v(\xi) ;\xi) $.  

\paragraph{$L^2$ optimality (Proper Generalized Decomposition methods).}
When one is interested in optimality in the norm $\Vert\cdot\Vert_2$, an $m$-dimensional subspace $V_m$ can be defined as the solution of the following optimization problem over the Grassmann manifold of subspaces of dimension $m$:
\begin{equation}
\min_{\dim(V_m)=m} \min_{v\in V_m\otimes L^2_\mu(\Xi)} \int_\Xi \Delta(v(\xi);\xi)^2d\mu(\xi),\label{optimVm}
\end{equation}
which can be equivalently written as 
 an optimization problem over the set of $m$-dimensional bases in $V$:
\begin{equation}
\min_{v_1,\hdots,v_m \in V} \min_{s_1,\hdots,s_m \in L^2_\mu(\Xi)} \int_\Xi \Delta(\sum_{i=1}^m v_i s_i(\xi);\xi)^2 d\mu(\xi).\label{optimbasisL2}
\end{equation}
This problem is an optimization problem over the set $\Rc_m$ of rank-$m$ tensors  which can be solved using optimization algorithms on low-rank manifolds such as an alternating minimization algorithm which consists in successively minimizing over the $v_i$ and over the $s_i$.
Assuming that $\Delta(\cdot;\xi)$ defines a distance to the solution $u(\xi)$ which  is uniformly equivalent to the one induced by the norm $\Vert \cdot\Vert_V$, i.e.
\begin{align}
\alpha_\Delta  \Vert u(\xi) - v \Vert_V \le \Delta(v;\xi) \le \beta_\Delta \Vert u(\xi) - v \Vert_V,\label{Delta_equivalence}
\end{align}
the resulting subspace $V_m$ is quasi-optimal in the sense that 
$$
\Vert u - P_{V_m} u \Vert_2 \le {c}_\Delta\min_{\dim(V_m)=m} \Vert u - P_{V_m} u \Vert_2 = {c}_\Delta d_m^{(2)}(u),
$$
with ${c}_\Delta=  \frac{\beta_\Delta}{\alpha_\Delta} $.
For linear problems where $ R(v ;\xi)  = A(\xi) v - b(\xi) \in W$ and $\Delta(v(\xi);\xi) = \Vert R(v(\xi) ;\xi) \Vert_W$, \eqref{Delta_equivalence} results from the property  \eqref{assumption_A} of the operator $A(\xi)$.
A suboptimal but constructive variant of algorithm \eqref{optimbasisL2} is defined by 
\begin{equation}
\min_{v_m \in V} \min_{s_1,\hdots,s_m \in L^2_\mu(\Xi)} \int_\Xi \Delta(\sum_{i=1}^m v_i s_i(\xi);\xi)^2 d\mu(\xi),\label{optimVm_lowrank}
\end{equation}
which is a greedy construction of the reduced basis $\{v_1,\hdots,v_m\}$. It yields a nested sequence of subspaces $V_m$.
This is one of the variants of Proper Generalized Decomposition methods (see \cite{NOU07,NOU08b,NOU10}).

Note that 
in practice, for solving \eqref{optimbasisL2} or \eqref{optimVm_lowrank} when $\Xi$ is not a finite set, one has either to approximate functions $s_i$ in a finite-dimensional subspace of $L^2_\mu(\Xi)$ \cite{NOU08b,NOU07}, or to approximate the integral over $\Xi$ by using a suitable quadrature \cite{Giraldi:2015:tobeornottobe}.

\paragraph{$L^\infty$ optimality (Reduced Basis methods).} 
When one is interested in optimality in the norm $\Vert\cdot\Vert_\infty$, a subspace $V_m$ could be defined by
$$
\min_{\dim(V_m)=m} \max_{\xi \in \Xi} \Delta(u_m(\xi);\xi),
$$
where $u_m(\xi)$ is some  projection of $u(\xi)$ onto $V_m$ 
(typically a Galerkin projection). A modification of the above definition consists in searching for 
spaces $V_m $ that are generated from evaluations of the solution at $m$ points selected in a subset $\Xi_K$ of $K$ points in $\Xi$ (a training set). In practice, this combinatorial optimization problem can be replaced by a greedy algorithm for the selection of points: given a set of interpolation points $\{\xi_*^1,\hdots,\xi_*^m\}$ and an approximation $u_m(\xi)$ in $V_m = \mathrm{span}\{u(\xi^1_*),\hdots,u(\xi^m_*)\}$, a new point $\xi^{m+1}_*$ is selected such that 
\begin{equation}
 \Delta(u_m(\xi^{m+1}_*);\xi^{m+1}_*) = \max_{\xi\in \Xi_K} \Delta(u_m(\xi);\xi).\label{greedy_selection_indicator}
\end{equation} 
This results in an adaptive interpolation algorithm
which was first introduced in \cite{PRU02}. It is the basic idea behind the so-called \emph{Reduced Basis methods} (see e.g. \cite{Patera:2007_book,Quarteroni:2011}). 
\\
Assuming that $\Delta$ satisfies \eqref{Delta_equivalence} and that the projection $u_m(\xi)$ verifies the quasi-optimality condition \eqref{quasioptimality}, the selection of interpolation points defined by \eqref{greedy_selection_indicator} is quasi-optimal in the sense that 
$$
\Vert u(\xi^{m+1}_*) - P_{V_m} u(\xi^{m+1}_*) \Vert_V \ge \gamma \max_{\xi\in \Xi_K} \Vert u(\xi) - P_{V_m} u(\xi) \Vert_V,
$$
with $\gamma =  c_\Delta^{-1}\inf_{\xi \in \Xi} c(\xi)^{-1}$. That makes \eqref{greedy_selection_indicator}  a suboptimal  version of the greedy algorithm 
 \eqref{eim} (so-called weak greedy algorithm). Convergence results for this algorithm can be found in \cite{Buffa:2012,binev:2011,DeVore:2013fk}, where the authors provide explicit comparisons  between the resulting error $\Vert u(\xi) - u_m(\xi) \Vert_V$ and the Kolmogorov $m$-width of $u(\Xi_K)$ which is the best achievable error  by projections onto $m$-dimensional spaces.  
 \\
The above definitions of interpolation points, and therefore of the resulting subspaces $V_m$, do not take into account explicitly the probability measure $\mu$. However, this measure is taken into account implicitly when working with a sample set  $\Xi_K$ drawn according the probability measure $\mu$. 
A construction that takes into account the measure explicitly has been proposed in \cite{Chen:2013}, where $\Delta(v;\xi)$ is replaced by 
$\omega(\xi) \Delta(v;\xi)$, with $\omega(\xi)$ a weight function depending on the probability measure $\mu$.
 
\section{Low-rank approximation of multivariate functions}

For the approximation of high-dimensional functions $u(\xi_1,\hdots,\xi_d)$, a standard approximation tool consists in searching for an expansion on a multidimensional basis obtained by tensorizing univariate bases:
\begin{equation}
u(\xi_1,\hdots,\xi_d) \approx \sum_{\nu_1=1}^{n_1} \hdots \sum_{\nu_d=1}^{n_d} a_{\nu_1,\hdots,\nu_d} \phi^1_{\nu_1}(\xi_1) \hdots \phi^d_{\nu_d}(\xi_d).    \label{tensorapprox}
\end{equation}
This results in an exponential growth of storage and computational complexities with the dimension $d$. Low-rank tensor methods aim at reducing the complexity by exploiting high-order low-rank structures of multivariate functions, considered as elements of tensor product spaces. This section presents basic notions about low-rank approximations of high-order tensors.
The reader is referred to the textbook \cite{HAC12} and the  surveys \cite{KOL09,Khoromskij:2012fk,Grasedyck:2013} for further details on the subject. 
For simplicity, we consider the case of a real-valued function $u:\Xi \to \Rbb$. The presentation naturally extends to the case $u:\Xi \to V$ by the addition of a new dimension.
\\
\par Here, we assume that $\Xi=\Xi_1\times \hdots \times \Xi_d$, and that $\mu$ is a product measure $\mu_1\otimes \hdots \otimes \mu_d$.  
Let $S_\nu$ denote a space of univariate functions defined on $\Xi_\nu$, $1\le \nu\le d$.
The elementary tensor product of functions $v^\nu \in S_\nu$ is defined 
by
$
(v^1\otimes \hdots \otimes v^d)(\xi_1,\hdots,\xi_d) = v^1(\xi_1)\hdots v^d(\xi_d).
$
The algebraic tensor space $S = S_1\otimes \hdots \otimes S_d$ is defined as the set of elements that can be written as a finite sum of elementary tensors, which means
\begin{align}
v(\xi_1,\hdots,\xi_d)   
= \sum_{i=1}^r a_i v_i^1 (\xi_1) \hdots v_i^d(\xi_d),
\label{canonicaldecomp}
\end{align}
for some $v_i^\nu \in S_\nu$, $a_i \in \Rbb$ and  $r\in \Nbb$. For the sake of simplicity, we consider that $S_\nu$ is a finite-dimensional approximation space in $L^2_{\mu_\nu}(\Xi_\nu)$
 (e.g. a space of polynomials, wavelets, splines...), with $\dim(S_\nu)=n_\nu \le n$, so that $S$ is a subspace of the algebraic tensor space $L^2_{\mu_1}(\Xi_1)\otimes \hdots \otimes L^2_{\mu_d}(\Xi_d)$ (whose completion is $L^2_\mu(\Xi)$).

\subsection{Tensor ranks and corresponding low-rank formats}
\begin{figure}
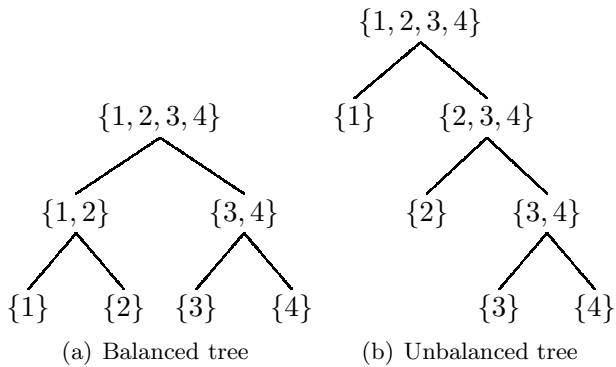

\centering
\subfigure[Balanced tree]{
\synttree[$\{1,2,3,4\}$[$\{1,2\}$[$\{1\}$][$\{2\}$]][$\{3,4\}$[$\{3\}$][$\{4\}$]]]\label{{figtreeHT}}}
\subfigure[Unbalanced tree]{
\synttree[$\{1,2,3,4\}$[$\{1\}$][$\{2,3,4\}$[$\{2\}$][$\{3,4\}$[$\{3\}$][$\{4\}$]]]]\label{figtreeTT}}
\caption{Examples of dimension partition trees over $D=\{1,\hdots,4\}$.}.
\label{figtrees}
\end{figure}

The \emph{canonical rank} of a tensor $v$ in $S$ is the minimal integer $m$ such that $v$ can be written under the form 
\eqref{canonicaldecomp}. An approximation of the form \eqref{canonicaldecomp} is called an approximation in \emph{canonical tensor format}. It has a storage complexity in $O(rnd)$. For order-two tensors, this is the standard and unique  notion of rank. 
For higher-order tensors, other notions of rank can be introduced, therefore yielding different types of rank-structured approximations. First, for a certain subset of dimensions $\alpha \subset D:=\{1,\hdots,d\}$ and its complementary subset $\alpha^c = D\setminus \alpha$, $S$ can be identified with the space $S_\alpha \otimes S_{\alpha^c}$ of order-two tensors, where $S_\alpha = \bigotimes_{\nu\in \alpha} S_\nu$. The $\alpha$-rank of a tensor $v$, denoted $\rank_\alpha(v)$, is then defined as the minimal integer $r_\alpha$ such that 
$$
v(\xi_1,\hdots,\xi_d) = \sum_{i=1}^{r_\alpha} v_i^{\alpha}(\xi_\alpha) v_i^{\alpha^c}(\xi_{\alpha^c}), 
$$
where $v_i^\alpha \in S_\alpha$ and $v_i^{\alpha^c} \in S_{\alpha^c}$ (here $\xi_\alpha$ denotes the collection of variables $\{\xi_\nu:\nu\in \alpha\}$), which is the standard and unique notion of rank for order-two tensors. Low-rank Tucker formats are then defined by imposing the $\alpha$-rank for a collection of subsets $\alpha \subset D$. The \emph{Tucker rank}  (or \emph{multilinear rank})  of a tensor is the tuple $ (\rank_{\{1\}}(v),\hdots,\rank_{\{d\}}(v)) $ in $\Nbb^d$. A tensor with Tucker rank bounded by $  (r_1,\hdots,r_d)$ can be written 
\begin{equation}
v(\xi_1,\hdots,\xi_d) = \sum_{i_1=1}^{r_1} \hdots \sum_{i_d=1}^{r_d} a_{i_1\hdots i_d} v_{i_1}^1(\xi_1)\hdots v_{i_d}^d(\xi_d),\label{tuckerformat}
\end{equation}
where $v_{i_\nu}^\nu \in S_\nu$ and where $a \in \Rbb^{r_1\times \hdots \times r_d}$ is a tensor of order $d$, called the \emph{core tensor}. An approximation of the form  \eqref{tuckerformat} is called an approximation in \emph{Tucker format}. It can be seen as an approximation in a tensor space $U_1\otimes \hdots\otimes U_d$, where $U_\nu=\mathrm{span}\{v_{i_\nu}^\nu\}_{i_\nu=1}^{r_\nu}$ is a $r_\nu$-dimensional subspace of $S_\nu$. The storage complexity of this format is in $O(rnd + r^d)$ (with $r = \max_\nu r_\nu$) and grows exponentially with the dimension $d$. Additional constraints on the ranks of $v$ (or of the core tensor $a$) have to be imposed in order to reduce this complexity. \emph{Tree-based (or Hierarchical) Tucker formats} \cite{HAC09,Falco2013} are based on a notion of rank associated with a dimension partition tree $T$ which is a tree-structured collection of subsets $\alpha$ in $D$, with $ D$ as the root of the tree and the  singletons $\{1\},\hdots,\{d\}$ as  the leaves of the tree (see Figure \ref{figtrees}). The \emph{tree-based (or Hierarchical) Tucker rank} associated with $T$ is then defined as the tuple $(\rank_\alpha(v))_{\alpha\in T}$. A particular case of interest is the \emph{Tensor Train (TT) format} \cite{OSE11} which is associated with a simple rooted tree $T$ with interior nodes $I = \{\{\nu,\hdots,d\}$, $1\le \nu \le d$\} (represented on Figure \ref{figtreeTT}). The \emph{$TT$-rank} of a tensor $v$ is defined as the tuple of ranks $(\rank_{\{\nu+1,\hdots,d\}}(v))_{1\le \nu\le d-1}$. A tensor $v$ with TT-rank bounded by $(r_1,\hdots,r_{d-1})$ can be written under the form
\begin{equation}
v(\xi_1,\hdots,\xi_d) = \sum_{i_1=1}^{r_1}\hdots \sum_{i_{d-1}=1}^{r_{d-1}} v_{1,i_1}^1(\xi_1)v_{i_1,i_2}^2(\xi_2)\hdots v_{i_{d-1},1}^d(\xi_d),
\label{TTformat}
\end{equation}
with $v_{i_{\nu-1},i_\nu}^\nu \in S_\nu$, which is very convenient in practice (for storage, evaluations and algebraic manipulations). The storage complexity for the TT format  is in 
$O(dr^2n)$ (with $r=\max_{\nu} r_\nu$).  

\subsection{Relation with other structured approximations}

Sparse tensor methods 
consist in searching for approximations of the form \eqref{tensorapprox} with only a few non-zero terms, that means approximations $v(\xi) = \sum_{i=1}^m a_i s_i(\xi)$ where the $m$ functions $s_i$ are selected (with either adaptive or non-adaptive methods) in the collection (\emph{dictionary}) of functions
$\Dc= \{ \phi_{k_1}^1(\xi_1) \hdots \phi_{k_d}^d(\xi_d): 1 \le k_\nu  \le n_\nu,  1\le \nu \le d\}$.  A typical choice consists in taking for $\Dc$ a basis of multivariate polynomials.  Recently, theoretical results have been obtained  on the convergence of sparse polynomial approximations of the solution of parameter-dependent equations (see \cite{2015arXiv150206797C}). Also, algorithms have been proposed that can achieve convergence rates comparable to the best $m$-term approximations.
For such a dictionary $\Dc$ containing rank-one functions, a sparse $m$-term approximation is a tensor with  canonical rank bounded by $m$ (usually lower than $m$, see \cite[Section 7.6.5]{HAC12}). Therefore, an
approximation in canonical tensor format  \eqref{canonicaldecomp} can be seen as a sparse $m$-term approximation where the $m$ functions are selected in the dictionary of all rank-one (separated) functions $\Rc_1 = \{s(\xi) = s^1(\xi_1)\hdots s^{d}(\xi_d) : s^\nu \in S^\nu\}$. 
Convergence results for best rank-$m$ approximations can therefore 
be deduced from convergence results for sparse tensor approximation methods. 

Let us mention some other standard structured approximations that are particular cases of low-rank tensor approximations.
First, a function   $v(\xi) = v(\xi_\nu)$ depending on a single variable $\xi_\nu$ is a rank-one (elementary) tensor. It has an $\alpha$-rank equal to $1$ for any subset $\alpha$ in $D$. A low-dimensional function $v(\xi) = v(\xi_\alpha)$ depending on a subset of variables $\xi_\alpha$, $\alpha \subset D$, has a $\beta$-rank equal to $1$ for any subset of dimensions $\beta$ containing $\alpha$ or such that $\beta\cap \alpha = \emptyset.$   An additive function $v(\xi) = v_1(\xi_1) + \hdots+ v_d(\xi_d)$, which is a sum of $d$ elementary tensors, is a tensor with canonical rank $d$. Also, such an additive function has $\rank_\alpha(v)\le2$ for any subset $\alpha\in D $, which means that it admits an exact representation in any Hierarchical Tucker format (including TT format) with a rank bounded by $(2,\hdots,2)$. 
 
\begin{rem}
Let us note that low-rank structures (as well as other types of structures) can be revealed only after a suitable change of variables. For example, let  $\eta = (\eta_1,\hdots,\eta_m)$ be the variables obtained by an affine transformation of variables $\xi$, with 
$\eta_i=\sum_{j=1}^d a_{ij}\xi_j+b_j$, $1\le j\le m$. Then the function 
$v(\xi) = \sum_{i=1}^m v_{i}(\eta_i) := \hat v(\eta)$, as a function $m$ variables, can be seen as a order-$m$ tensor with canonical rank less than $m$. This type of approximation corresponds to the projection pursuit regression model.
 \end{rem}

\subsection{Properties of low-rank formats}

Low-rank tensor approximation methods consist in searching for approximations in a  subset of tensors 
$$
\Mc_{\le r} = \{v : \rank(v) \le r\},
$$
where different notions of rank yield different approximation formats. A first important question is to characterize the approximation power of low-rank formats, which means to quantify  the best approximation error 
$$
\inf_{v \in \Mc_{\le r}} \Vert u - v \Vert := \sigma_r
$$
for a given class of functions $u$ and a given low-rank format. A few convergence results have been obtained for functions with standard Sobolev regularity (see e.g. \cite{Schneider201456}). 
An open and challenging problem is to characterize  approximation classes of the different low-rank formats, which means the class of functions 
$u$ such that $\sigma_r$ has a certain (e.g. algebraic or exponential) decay with $r$.

The characterization of topological and geometrical properties of subsets $\Mc_{\le r}$ is important for different purposes such as proving the existence of best approximations in $\Mc_{\le r}$ or  deriving algorithms.
For $d\ge 3$, the subsets $\Mc_{\le r}$ associated with the notion of canonical rank are not closed and  best approximation problems in $\Mc_{\le r}$ are ill-posed.
Subsets of tensors associated with the notions of tree-based (hierarchical) Tucker ranks have better properties. Indeed, they are closed sets, which ensures the existence of best approximations. Also, they are differentiable manifolds \cite{2015arXiv150503027F,uschmajew2013geometry,Holtz:2012fk}. This has useful consequences for optimization \cite{uschmajew2013geometry} or for the projection of dynamical systems on these manifolds \cite{Lubich:2013}.  

For the different notions of rank introduced above, an interesting property of the corresponding low-rank subsets is that they admit simple parametrizations
\begin{equation}
\Mc_{\le r} = \{v = F(p_1,\hdots,p_\ell) : p_k \in P_k\},\label{subsetparam}
\end{equation}
where $F$ is a multilinear map and $P_k$ are vector spaces, and where the number of parameters $\ell$ is in $O(d)$. An optimization problem on $\Mc_{\le r}$ can then be reformulated as an optimization problem on $P_1\times \hdots \times P_\ell$, for which simple alternating minimization algorithms (block coordinate descent) can be used \cite{2015arXiv150600062E}.

\section{Low-rank approximation from samples of the function}
Here we present some strategies for the construction of low-rank approximations of a multivariate function $u(\xi)$ from point evaluations of the function.

\subsection{Least-squares}
Let us assume that $u\in L^2_\mu(\Xi)$. Given 
a set of $K$ samples $\xi^1,\hdots,\xi^K$ of $\xi$ drawn according the probability measure $\mu$, a least-squares approximation of $u$ in a low-rank subset $\Mc_{\le r}$ is defined by 
\begin{equation}
\min_{v \in \Mc_{\le r}} \frac{1}{K} \sum_{k=1}^K (u(\xi^k) - v(\xi^k))^2.\label{least-squares}
\end{equation}
Using a multilinear parametrization of low-rank tensor subsets (see \eqref{subsetparam}), the optimization problem \eqref{least-squares} can be solved using alternating minimization algorithms, each iteration corresponding to a standard least-squares minimization for  linear approximation \cite{BEL11,Doostan:2013,2013arXiv1305.0030C}. An  open question concerns the analysis of the number of samples which is required for a stable approximation in a given low-rank format. Also, standard regularizations can be introduced, such as sparsity-inducing regularizations \cite{2013arXiv1305.0030C}. In this statistical framework, cross-validation methods can  be used for the selection of tensor formats, ranks and approximation spaces $S_k$ (see \cite{2013arXiv1305.0030C}). 

\begin{rem}
Note that for other objectives in statistical learning (e.g. classification), \eqref{least-squares} can be replaced by  
\begin{equation*}
\min_{v \in \Mc_{\le r}} \frac{1}{K} \sum_{k=1}^K \ell(u(\xi^k),v(\xi^k)),
\end{equation*}
where $\ell$ is a so-called \emph{loss function} measuring a certain distance between $u(\xi^k)$ and the approximation $v(\xi^k)$.
\end{rem}

\subsection{Interpolation/Projection}

Here we present interpolation and projection methods for the approximation of $u$ in $S=S_1\otimes \hdots \otimes S_d$. Let $\{\phi^\nu_{k_\nu}\}_{k_\nu \in \Lambda_{n_\nu}}$ be a basis of $S_\nu$, with $\Lambda_{n_\nu}= \{1,\hdots,n_\nu\}$.  
If $\{\phi^\nu_{k_\nu}\}_{k_\nu \in \Lambda_{n_\nu}}$ is a set of interpolation functions associated with a set of points  $\{\xi_\nu^{k_\nu}\}_{k_\nu\in \Lambda_{n_\nu}}$ in $\Xi_\nu$,  then $\{\phi_k(\xi) = \phi_{k_1}^1(\xi_1)\hdots \phi_{k_d}^d(\xi_d)\}_{k\in \Lambda}$, $\Lambda=\Lambda_{n_1}\times \hdots \times \Lambda_{n_d}$, is a set of interpolation functions associated with 
the tensorized grid $\{\xi^k=(\xi_1^{k_1},\hdots,\xi_d^{k_d})\}_{k\in \Lambda}$ composed of $N=\prod_{\nu=1}^d n_\nu$ points. An interpolation $u_N$ of $u$ is then given by 
 $$
u_N(\xi) = \sum_{k\in \Lambda} u(\xi^k) \phi_k(\xi), 
 $$
 so that $u_N$ is completely characterized by the order-$d$ tensor $U \in \Rbb^{n_1}\otimes \hdots \otimes \Rbb^{n_d}$ whose components $U_{k_1,\hdots,k_d}  = u(\xi_{1}^{k_1},\hdots,\xi_d^{k_d})$ are the evaluations of $u$ on the interpolation grid. Now, if $\{\phi^\nu_{k_\nu}\}_{k_\nu \in \Lambda_{n_\nu}}$  is an orthonormal basis of $S_\nu$ (e.g. orthonormal polynomials) and if $\{(\xi_\nu^{k_\nu},\omega^\nu_{k_\nu})\}_{k_\nu \in \Lambda_{n_\nu}}$ is a quadrature rule on $\Xi_\nu$ (associated with the measure $\mu_\nu$), an approximate $L^2$-projection of $u$ can also be defined by 
\begin{align*}
u_N(\xi) = \sum_{k \in \Lambda} u_{k} \phi_k(\xi), \quad 
u_k = \sum_{k\in \Lambda} \omega_k u(\xi^k)\phi_k(\xi^k), 
\end{align*}
with $\omega_k = \omega^1_{k_1}\hdots \omega_{k_d}^d$.
Here again, $u_N$ is completely characterized by the order-$d$ tensor $U$ whose components are the evaluations of $u$ on the $d$-dimensional quadrature grid.

Then, low-rank approximation methods can be used in order to obtain an approximation of $U$ using only a few entries of the tensor (i.e. a few evaluations of the function $u$). 
This is related to the problem of tensor completion. A possible approach consists in evaluating some entries of the tensor taken at random and then in reconstructing the tensor by the minimization of a least-squares functional (this is the algebraic version of the least-squares approach described in the previous section) or dual approaches using regularizations of rank minimization problems (see \cite{2014arXiv1404.3905R}). In this statistical framework, a challenging question is to determine the number of samples required for a stable reconstruction of low-rank approximations in different tensor formats (see  \cite{2014arXiv1404.3905R} for first results). An algorithm has been introduced in \cite{Espig:2009yf} for the approximation in canonical format, using least-squares minimization with a structured set of entries selected adaptively.  
Algorithms have also been proposed for an adaptive construction of low-rank approximations of $U$ in  Tensor Train format \cite{OSE10} or  Hierarchical Tucker format \cite{Ballani2013639}.   These algorithms are extensions of Adaptive Cross Approximation (ACA) to high-order tensors and provide approximations that interpolate the tensor $U$ at some adaptively chosen entries.

\section{Low-rank tensor methods for parameter-dependent equations}
Here,  we present numerical methods for the direct computation of low-rank approximations of the solution $u$ of a parameter-dependent equation $R(u(\xi);\xi)=0$, where $u$ is seen as a two-order tensor in $V\otimes L^2_\mu(\Xi)$, or as a higher-order tensor by exploiting an additional tensor space structure of $ L^2_\mu(\Xi)$ (for a product measure $\mu$). 

\begin{rem}
When exploiting only the order-two tensor structure, the methods presented here are closely related to projection-based model reduction methods. Although they provide a directly exploitable low-rank approximation of $u$, 
they can also be used for the construction of a low-dimensional subspace in $V$ (a candidate for projection-based model reduction) which is extracted from the obtained low-rank approximation. 
\end{rem} 
 
\subsection{Tensor-structured equations}
Here, we describe how the initial equation can be reformulated as a tensor-structured equation.
In practice, a  preliminary discretization of functions defined on $\Xi$ is required.  
A possible discretization consists in introducing a $N$-dimensional approximation space $S$ in 
$L^2_\mu(\Xi)$ (e.g. a polynomial space) and a standard Galerkin projection of the solution onto $V\otimes S$ (see e.g. \cite{MAT05,NOU10}). The resulting approximation can then be identified with a two-order tensor $\mathbf{u}$ in $V\otimes \Rbb^N$. When $S$ is the tensor product of $n_\nu$-dimensional spaces $S_\nu$, $1\le \nu\le d$, the resulting approximation can be identified with a higher-order tensor $\mathbf{u}$ in $V\otimes \Rbb^{n_1}\otimes \hdots \otimes \Rbb^{n_d}$.
Another simple discretization consists in considering only a finite (possibly large) set of $N$ points in $\Xi$ (e.g. an interpolation grid) and the corresponding finite set of equations $R(u(\xi^k);\xi^k)  = 0, \; 1\le k\le N,$ on the set of samples of the solution $ (u(\xi^k))_{k=1}^N \in V^N$, which can be identified with a tensor 
 $\mathbf{u} \in V\otimes \Rbb^N$. If the set of points is obtained by the tensorization of unidimensional grids with $n_\nu$ points, $1\le \nu\le d$, then $(u(\xi^k))_{k=1}^N$ can be identified with a higher-order tensor $\mathbf{u}$ in $V\otimes \Rbb^{n_1}\otimes \hdots \otimes \Rbb^{n_d}$. 
Both types of discretization yield an equation
\begin{equation}
\mathbf{R}(\mathbf{u})=0,\label{residualeq}
\end{equation} 
with 
$$\mathbf{u} \text{ in } V\otimes \Rbb^N \text{ or } V\otimes \Rbb^{n_1}\otimes \hdots \otimes \Rbb^{n_d}.$$
\\
In order to clarify the structure of equation \eqref{residualeq}, let us consider the case 
of a linear problem where $R(v;\xi) = A(\xi) v - b(\xi)$ and assume that $A(\xi)$ and $b(\xi)$ have low-rank (or affine) representations of the form 
\begin{equation}
A(\xi) = \sum_{i=1}^L {A_i} \alpha_i(\xi) \quad \text{and} \quad   b(\xi)=  \sum_{i=1}^R b_i  \beta_i(\xi).\label{eq:affine_A_b_tensor}
\end{equation}
Then \eqref{residualeq} takes the form of a tensor-structured equation $ \mathbf{A}\mathbf{u}-\mathbf{b}=0,$ with 
\begin{align}
 \mathbf{A}= \sum_{i=1}^L A_i \otimes \tilde A_i \quad \text{and} \quad  \mathbf{b}= \sum_{i=1}^R b_i \otimes \tilde b_i,\label{tensorstructuredAb}
\end{align}
where for the second type of discretization, $\tilde A_i \in \Rbb^{N\times N}$ is a diagonal matrix whose diagonal is the vector of evaluations of $\alpha_i(\xi)$ at the sample points, and $\tilde b_i \in \Rbb^N$ is the vector of evaluations of $\beta_i(\xi)$ at the sample points. If $A(\xi)$ and $b(\xi)$ have higher-order low-rank representations of the form 
\begin{align}
A(\xi) = \sum_{i=1}^L A_i \alpha_{i}^1(\xi_1)\hdots \alpha_{i}^d(\xi_d)\quad \text{and} \quad  b(\xi)=  \sum_{i=1}^R b_i  \beta_i^1(\xi_1)\hdots \beta_i^d(\xi_d),  \label{eq:affine_A_b_high-order}
\end{align} 
then  \eqref{residualeq} takes the form of a tensor-structured equation $ \mathbf{A}\mathbf{u}-\mathbf{b}=0$ on $\mathbf{u}$ in $V\otimes \Rbb^{n_1}\otimes \hdots \otimes \Rbb^{n_d}$, with
\begin{align}
 \mathbf{A}= \sum_{i=1}^L A_i \otimes \tilde A_i^1 \otimes \hdots \otimes \tilde A_i^d \quad \text{and} \quad  \mathbf{b}= \sum_{i=1}^R b_i \otimes \tilde b_i^1 \otimes \hdots \otimes \tilde b_i^d,\label{tensorstructuredAb-highorder}
\end{align}
where for the second type of discretization (with a tensorized grid in $\Xi_1\times \hdots \times \Xi_d$), $\tilde A_i^\nu \in \Rbb^{n_\nu\times n_\nu}$ is a diagonal matrix whose diagonal is the vector of evaluations of $\alpha_i^\nu(\xi_\nu)$ on the unidimensional grid in $\Xi_\nu$, and $\tilde b_i^\nu \in \Rbb^{n_\nu}$ is the vector of evaluations of $\beta_i^\nu(\xi_\nu)$ on this grid (for the first type of discretization, see  \cite{NOU10} for the definition of tensors $\mathbf{A}$ and $\mathbf{b}$).
Note that when $A(\xi)$ and $b(\xi)$ do not have low-rank representations \eqref{eq:affine_A_b_tensor} or  \eqref{eq:affine_A_b_high-order} (or any other higher-order low-rank representation), then a preliminary approximation step is required in order to obtain such approximate representations (see e.g. \cite{BAR02,Casenave:2014fk,2015arXiv150303210D}). This is crucial for reducing the computational and storage complexities.

\subsection{Iterative solvers and low-rank truncations}
A first solution strategy consists in using standard iterative solvers (e.g. Richardson, conjugate gradient, Newton...) with efficient low-rank truncation methods of the iterates \cite{Ballani:2013vn,Kressner:2011ys,KHO11,MAT12,2015arXiv150107714B,Bachmayr:2015fk}. 
A simple iterative algorithm takes the form 
$\mathbf{u}^{k+1} = M(\mathbf{u}^k),$ where $M$ is an iteration map involving simple algebraic operations between tensors (additions, multiplications) which requires the implementation of a  tensor algebra. 
Low-rank truncation methods can be systematically used for limiting the storage complexity and the computational complexity of algebraic operations. This results in approximate iterations 
$
\mathbf{u}^{k+1} \approx M(\mathbf{u}^k)
$
and the resulting algorithm 
 can be analyzed as an inexact version (or perturbation) of the initial algorithm (see e.g. \cite{HAC08}).

As an example, let us consider a linear tensor-structured problem $\mathbf{A}\mathbf{u}-\mathbf{b}=0$. An approximate Richardson algorithm takes the form 
$$\mathbf{u}^{k+1} = \Pi_\epsilon(\mathbf{u}^k + \alpha (\mathbf{b}-\mathbf{A} \mathbf{u}^k)),$$ where $\Pi_\epsilon$ is a map which associates to a tensor $\mathbf{w}$ a low-rank approximation $\Pi_\epsilon(\mathbf{w})$ such that 
$\Vert \mathbf{w} - \Pi_\epsilon (\mathbf{w})\Vert_{p} \le \epsilon \Vert \mathbf{w} \Vert_p$, with $p=2$ or $p=\infty$ depending on the desired control of the error (mean-square or uniform error control over the parameter set). Provided some standard assumptions on the operator $\mathbf{A}$ and the parameter $\alpha$, the generated sequence $\mathbf{u}^k$ is such that
$\lim\sup_{k\to\infty} \Vert \mathbf{u} - \mathbf{u}^k \Vert_p \le C(\epsilon) $ with $C(\epsilon) \to 0$ as $\epsilon\to 0$.  
For $p=2$, efficient low-rank truncations of a tensor can be obtained using SVD for an order-two tensor or generalizations of SVD for higher-order tensor formats \cite{OSE09,GRA10}. A selection of the ranks based on the singular values of (matricisations of) the tensor allows a control of the error.  In  \cite{2015arXiv150107714B}, the authors propose an alternative truncation strategy based on soft thresholding. For $p=\infty$, truncations can be obtained with an \green{\emph{Adaptive Cross Approximation}} algorithm (or Empirical Interpolation Method) \cite{Bebendorf:2014uq} for an order-two tensor or with extensions of this algorithm for higher-order tensors \cite{OSE10}. Note that low-rank representations of the form  \eqref{tensorstructuredAb} or \eqref{tensorstructuredAb-highorder} for $\mathbf{A}$ and $\mathbf{b}$ are crucial since they ensure that algebraic operations between tensors can be done with a reduced complexity. 
Also, iterative methods usually require good preconditioners. In order to maintain a low computational complexity, these preconditioners must also admit  low-rank representations.

\subsection{Optimization on low-rank manifolds}
Another solution strategy consists in directly computing a low-rank approximation by minimizing some functional $\Jc$ whose minimizer on $V\otimes \Rbb^N$ (or $V\otimes \Rbb^{n_1}\otimes \hdots \otimes \Rbb^{n_d}$) is the solution of equation \eqref{residualeq}, i.e. by solving 
\begin{equation}
\min_{\mathbf{v}\in \Mc_{\le r}}\Jc(\mathbf{v}),\label{minJRr}
\end{equation}
where $\Mc_{\le r}$ is a low-rank manifold. 
There is a natural choice of functional for problems where 
\eqref{residualeq} corresponds to the stationary condition of a functional $\Jc$ \cite{Giraldi:2015:tobeornottobe}.  Also, $\Jc(\mathbf{v})$ can be taken as a certain norm of the residual $\mathbf{R}(\mathbf{v})$.
 For linear problems, choosing $$\Jc(\mathbf{v}) = \Vert \mathbf{A} \mathbf{v} - \mathbf{b} \Vert^2$$  yields a quadratic optimization problem over a low-rank manifold. Optimization problems on low-rank manifolds  
 can be solved either by using algorithms which exploit the manifold structure (e.g. Riemannian optimization) or by using simple alternating minimization algorithms given a parametrization \eqref{subsetparam} of the low-rank manifold. Under the assumption that $\Jc$ satisfies
$
\alpha \Vert \mathbf{u}-\mathbf{v} \Vert_2 \le \Jc(\mathbf{v}) \le \beta \Vert \mathbf{u}-\mathbf{v} \Vert_2, 
$
the solution $\mathbf{u}_r$ of \eqref{minJRr} is quasi-optimal in the sense that 
$$
\Vert \mathbf{u}-\mathbf{u}_r\Vert_2 \le \frac{\beta}{\alpha} \min_{\mathbf{v} \in \Mc_{\le r}} \Vert \mathbf{u}-\mathbf{v} \Vert_2,
$$
where $\beta/\alpha$ is the condition number of the operator $A$. Here again, the use of preconditioners allows us to better exploit the approximation power of a given low-rank manifold $\Mc_{\le r}$.

Constructive algorithms are also possible, the most prominent algorithm being a greedy algorithm which consists in computing a sequence of low-rank approximations $u_m$ obtained by successive rank-one corrections, i.e.
$$
\Jc(\mathbf{u}_m) = \min_{w\in \Rc_1} \Jc(\mathbf{u}_{m-1}+\mathbf{w}).
$$
This algorithm is a standard greedy algorithm \cite{2012arXiv1206.0392T} with a dictionary of rank-one (elementary)  tensors  $\Rc_1$, which was first introduced in \cite{NOU07} for the solution of parameter-dependent (stochastic) equations.
 Its 
convergence has been established under standard assumptions for convex optimization problems \cite{CAN11,Falco:2012fk}. 
The utility of this algorithm is that it is adaptive and it only requires the solution of optimization problems on the low-dimensional manifold $\Rc_1$. However, for many practical problems, greedy constructions of low-rank approximations in canonical low-rank format are observed to converge slowly. 
Improved constructive algorithms which better exploit the tensor  structure of the problem have been proposed \cite{NOU10} and convergence results are also available for some general convex optimization problems \cite{Falco:2012fk}.

\begin{rem} 
The above algorithms can also be used within iterative algorithms for which an iteration takes the form 
$$
\mathbf{C}_k \mathbf{u}^{k+1} = \mathbf{F}_k(\mathbf{u}^k),
$$
where $\mathbf{F}_k(\mathbf{u}^k)$ can be computed with a low complexity using low-rank tensor algebra (with potential low-rank truncations), but where the inverse 
of the operator $\mathbf{C}_k$ is not known explicitly, so that 
$\mathbf{C}_k^{-1}\mathbf{F}_k(\mathbf{u}^k)$ cannot be obtained from simple algebraic operations. Here, a low-rank approximation of $\mathbf{u}^{k+1}$ can be computed using the above algorithms with the residual-based functional 
$\Jc(\mathbf{v}) = \Vert \mathbf{C}_k \mathbf{v} - \mathbf{F}_k(\mathbf{u}^k) \Vert^2$.
\end{rem}

\section*{Concluding remarks}

Low-rank tensor methods have emerged as a very powerful tool for the solution of high-dimensional problems arising in many contexts, and in particular in uncertainty quantification. However, there remain many challenging issues to address for a better understanding of this type of approximation and for a  diffusion of these methods in a wide class of applications. 
From a theoretical point of view, open questions include 
the characterization of the approximation classes of a given low-rank tensor format, which means the class of functions for which a certain type of convergence (e.g. algebraic or exponential) can be expected, and also the characterization of 
the problems yielding solutions in these approximation classes. Also, quantitative results  on the approximation of a low-rank function (or tensor) from samples of this function (or entries of this tensor) could allow us to answer some practical issues such as the determination of the  number of samples required for a stable approximation in a given low-rank format or the design of sampling strategies which are adapted to particular low-rank  formats.  
From a numerical point of view, challenging issues include the development of efficient algorithms for global optimization on low-rank manifolds, with guaranteed convergence properties,  and the development of adaptive algorithms for the construction of controlled low-rank approximations, with an adaptive selection of ranks and potentially of the tensor formats (e.g. based on tree optimization for tree-based formats).
From a practical point of view, low-rank tensor methods exploiting model equations (Galerkin-type methods) are often seen as ``intrusive methods'' in the sense that they require (a priori) the development of specific softwares.  An important issue is then to develop weakly intrusive implementations of these methods which may allow  the use of existing computational frameworks and which would therefore contribute to a large diffusion of these methods.  


\end{document}